\documentclass[12pt]{article}
\usepackage[english,russian]{babel}
\usepackage{amsmath, amscd, amsthm, a4,  amssymb}
\usepackage[dvips]{graphicx}
\begin{document}

УДК  517.956

\begin{center}

{\bf О КОРРЕКТНОЙ РАЗРЕШИМОСТИ ЗАДАЧИ ДИРИХЛЕ ДЛЯ ОБОБЩЕННОГО УРАВНЕНИЯ МАНЖЕРОНА С НЕГЛАДКИМИ КОЭФФИЦИЕНТАМИ

\

И.Г.Мамедов}

Институт Кибернетики им А.И.Гусейнова НАН Азербайджана, AZ 1141, Азербайджан, Баку, ул. Б.Вагабзаде, 9

Email: ilgar-mammadov@rambler.ru

\end{center}

\begin{abstract}

{\it В работе для одного псевдопараболического  уравнения четвертого  порядка с негладкими коэффициентами рассмотрена задача Дирихле с неклассическими условиями, не требующими условий согласования.
Обоснована эквивалентность этих условий  классическим краевым условием в случае, если решение поставленной задачи ищется в изотропном пространстве С.Л.Соболева.
Решение осуществляется редукцией к системе уравнений Фредгольма, корректная разрешимость которых устанавливается при негладких условиях на коэффициенты уравнения на основе метода интегральных представлений.}\\

{\bf Ключевые слова:} задача Дирихле, псевдопараболические уравнения, уравнения с негладкими коэффициентами.
\end{abstract}

{\bf Введение}

Первая краевая задача или задача Дирихле, (т.е. задача в которой носителем данных является замкнутый контур) хорошо известная для дифференциальных уравнений эллиптического типа, является одной из основных краевых  задач математической физики [1]. С этой точки зрения, эта работа посвящена актуальным проблемам  математической физики.\\

{\bf 1. Постановка задачи.}
Рассмотрим обобщенное уравнение Манжерона
\begin{equation} \label{GrindEQ___A201_A20_}
{\left(V_{2, 2}  u\right) \left(x, y\right) \equiv  D_{x}^{2}  D_{y}^{2}  u \left(x , y\right) + a _{2 , 1}  \left(x , y\right) D_{x}^{2}  D_{y} u \left(x , y\right) + }
$$$$
 {+a_{1, 2 } \left(x , y\right) D_{x}  D_{y}^{2}  u \left(x , y\right) +a _{2 , 0}  \left(x , y\right) D_{x}^{2}  u \left(x , y\right) +}
$$$$
  {+a _{0,2 } \left(x , y\right) D_{y}^{2}  u \left(x , y\right) + \sum _{i=0}^{1}  \sum _{j=0}^{1}  a_{i , j } (x,y)D_{x}^{i}  D_{y}^{j}  u \left(x , y\right) = Z_{2 , 2 } \left(x , y\right) \in L_{p}  \left(G\right)}
\end{equation}

Здесь $u \left(x,y\right)$-искомая функция определенная на $G$;  $a_{i , j } $$\left(x , y\right)$ заданные измеримые функции на $G= G_{1} \times  G_{2} ,$ где $G_{j} =\left(0,h_{j}\right),$ $j=1,2,$ $Z_{2, 2} \left(x, y\right)$ заданная измеримая функция на $G;$ $D_{t}  \equiv  \partial /\partial $ $t-$оператор обобщенного дифференцирования в смысле С.Л. Соболева, $D_{t}^{0} $-оператор тождественного преобразования.

Уравнение (1) является гиперболическим уравнением, которое обладает двумя действительными характеристиками $x=const, y=const,$ первая и вторая из которых двукратная. Обобщенное уравнение Манжерона -- это одно из основных дифференциальных уравнений математической биологии [2]. Уравнения вида (1) в работе А.П.Солдатова и М.Х.Шханукова [3] названы псевдопараболическими. Рассматриваемое уравнение является обобщением многих модельных уравнений некоторых процессов (обобщенное уравнение влагопереноса, телеграфное уравнение, уравнение колебания струны, уравнение теплопроводности, уравнение Аллера и.т.д). Кроме того, это уравнение являющееся обобщением уравнения Буссинеска -- Лява [4], описывающего продольные волны в тонком упругом стержне с учетом эффектов поперечной инерции.

В работе рассматривается уравнение (1) в общем случае, когда коэффициенты  $a_{i , j } $$\left(x, y\right)$ являются негладкими функциями, удовлетворяющими лишь следующим условиям:
\[
\begin{array}{l} {a_{2, j} \left(x, y\right) \in  L_{\infty }^{x} ,_{p}^{y} \left(G\right),\quad j=\overline{0.1} ;} \\ {a_{i, 2}  \left(x, y\right) \in  L_{p}^{x} ,_{\infty }^{y} \left(G\right), \quad i=\overline{0.1 };} \\ {a_{i, j} \left(x, y\right) \in  L_{p} \left(G\right),\quad i={\overline{0, 1};}\quad   j=\overline{0.1} .} \end{array}
\]

При этом важным принципиальным моментом является то, что рассматриваемое уравнение обладает негладкими коэффициентами, которые удовлетворяют только некоторым условиям типа $p-$интегрируемости и ограниченности т.е. рассмотренный псевдопараболический дифференциальный оператор $V_{2, 2} $ не имеет традиционного сопряженного оператора. Иначе говоря, функция Римана для такого  уравнения не может быть исследована классическим методом характеристик.

При этих условиях решение $u(x, y)$  уравнения (1) будем искать в изотропном пространстве С.Л.Соболева
\[
W_{p}^{\left(2, 2\right)}  \left(G\right) \equiv  \left\{ u \left(x, y\right) : D_{x}^{i}  D_{y}^{j}  u \left(x, y\right) \in  L_{p}  \left(G\right), \quad i=\overline{0, 2},\quad j=\overline{0, 2}\right\},
\]
где $1\le p\le \infty $. Норму в пространстве С.Л.Соболева $W_{p}^{\left(2, 2\right)}  \left(G\right)$ будем определять равенством:
\[
\left\|  u \left(x, y\right) \right\|_{ W_{p}^{\left(2, 2\right)}  \left(G\right)} = \sum _{i=0}^{2}  \sum _{j=0}^{2}  \left\|  D_{x}^{i}  D_{y}^{j}  u \left(x, y\right) \right\|_{L_{p } \left(G\right)}.
\]

Для уравнения (1) условия Дирихле классического вида [5] можно задать в виде
\begin{equation} \label{GrindEQ__2_}
\left\{\begin{array}{l} {u \left(0, y\right) = \varphi _{1}  \left(y\right),   u \left(x, 0\right) = \psi _{1}  \left(x\right),} \\ {u \left(h_{1} , y\right) = \varphi _{2}  \left(y\right),   u \left(x, h_{2} \right) = \psi _{2}  \left(x\right),} \end{array}\right.
\end{equation}
где $\varphi _{j}  \left(y\right)$ и $\psi _{j } \left(x\right), j=\overline{1,2}$ - заданные измеримые функции на $G$. Очевидно, что случае условий \eqref{GrindEQ__2_} заданные функции кроме условий
\[
\varphi _{j}  \left(y\right) \in  W_{p}^{\left(2\right)}  \left(G_{2} \right) \equiv  \left\{\widetilde{\varphi } \left(y\right) : D_{y}^{j}  \widetilde{\varphi } \left(y\right)\in  L_{p}  \left(G_{2} \right),\ \ j=\overline{0, 2}\right\},\ \ 1\le p\le \infty ,\ \ j=1, 2;
\]
\[
\psi _{j}  (x) \in  W_{p}^{\left(2\right)}  \left(G_{1} \right) \equiv  \left\{ \widetilde{\psi } \left(x\right) : D_{x}^{i}  \widetilde{\psi } \left(x\right) \in  L_{p}  \left(G_1\right), i=\overline{0, 2}\right\} ,\ \  1\le p\le \infty ,\ \ j=\overline{1, 2}.
\]
должны удовлетворять также следующим условиям согласования:
\begin{equation} \label{GrindEQ__3_}
\left\{\begin{array}{l} {\varphi _{1}  \left(0\right)=\psi _{1}  \left(0\right),         \varphi _{2}  \left(h_{2} \right) = \psi _{2}  \left(h_{1} \right) ,} \\ {\varphi _{1}  \left(h_{2} \right) = \psi _{2}  \left(0\right),        \varphi _{2}  \left(0\right) = \psi _{1}  \left(h_{1} \right)}. \end{array}
\right.
\end{equation}

Рассмотрим следующие неклассические краевые условия:
\begin{equation} \label{GrindEQ__4_}
\left\{\begin{array}{l}
{V_{0, 0}  u \equiv  u \left(0, 0\right) = Z_{0, 0}  \in \mathbb{R} } \\
{V_{1, 0}  u \equiv  u_{x}  \left(0, 0\right) = Z_{1, 0}  \in \mathbb{R}} \\
{V_{0, 1}  u \equiv  u_{y}  \left(0, 0\right) = Z_{0, 1}  \in \mathbb{R}} \\
{\left(V_{2, 0}  u\right) \left(x\right) \equiv  u_{xx } \left(x, 0\right) =Z_{2, 0} \left(x\right) \in  L_{p}  \left(G_{1} \right) ;  } \\
{\left(V_{0, 2}  u\right) \left(y\right) \equiv  u_{yy } \left(0, y\right) =Z_{0, 2} \left(y\right) \in  L_{p}  \left(G_{2} \right) ;} \\
{V_{0, 0}^{\left(h_{ 1} \right)}  u \equiv  u \left(h_{1} , 0\right) = Z_{0, 0}^{\left(h_{ 1} \right)}  \in \mathbb{R} } \\ {V_{0, 1}^{\left(h_{ 1} \right)}  u \equiv  u_{y}  \left(h_{1} , 0\right) = Z_{0, 1}^{\left(h_{ 1} \right)}  \in \mathbb{R}  ;} \\
{V_{0 , 2}^{\left(h_{ 1} \right)}  u \equiv  u_{yy }  \left(h_{1, } y\right) = Z_{0, 2}^{\left(h_{ 1} \right)}  \left(y\right) \in  L_{p }  \left(G_{2} \right);} \\ {V_{0 , 0}^{\left(h_{ 2} \right)}  u \equiv   u \left(0_{, } h_{2} \right) = Z_{0, 0}^{\left(h_{ 2} \right)} \in \mathbb{R}  ; } \\
{V_{1, 0}^{\left(h_{ 2} \right)}  u\equiv  u_{x}  \left(0, h_{2} \right)= Z_{1, 0}^{(h_{ 2} ) }  \in \mathbb{R}  ;}
\\ { \left(V_{2, 0}^{\left(h_{ 2} \right)}  u\right) \left(x\right) \equiv  u_{xx }  \left(x, h_{2} \right) = Z_{2, 0}^{\left(h_{ 2} \right)}  \left(x\right) \in  L_{p}  \left(G_{1} \right)}. \end{array}\right.
\end{equation}

Если функция $u  \in   W_{p}^{(2, 2)}  \left(G\right)$ является решением задачи Дирихле классического вида (1), \eqref{GrindEQ__2_}, то она является также решением задачи (1), \eqref{GrindEQ__4_} для $Z_{i, j} $ и $Z_{i, j}^{\left(h_k\right)} ,$ определяемых следующими равенствами:
\[
\begin{array}{c} {Z_{0, 0}  = \varphi _{1}  \left(0\right) = \psi _{1}  \left(0\right) ; \ \  Z_{1, 0}  = \psi _{1}^{'}  \left(0\right) ;}\ \  {Z_{0, 1}  = \varphi _{1}^{'}  \left(0\right) ; \ \   Z_{2, 0} (x) =  \psi _{1}^{''}  \left(x\right) ;} \\ Z_{0, 2}  \left(y\right) = \varphi _{1}^{''}  \left(y\right) ; \ \ \   Z_{0, 0}^{\left(h_{ 1} \right)}  = \varphi _{2}  \left(0\right) = \psi _{1}  \left(h_{1} \right) ; \ \ {Z_{0, 1}^{\left(h_{ 1} \right)}  = \varphi _{2}^{'}  \left(0\right) ;  }\\  Z_{0, 0}^{\left(h_{ 2} \right)}  = \psi _{2}  \left(0\right) = \varphi _{1}  \left(h_{2} \right) ; \ {Z_{1, 0}^{\left(h_{ 2} \right)}  = \psi _{2}^{'}  \left(0\right) ; \  Z_{2,0}^{\left(h_{2} \right)}  \left(x\right) = \psi _{2}^{''}  \left(x\right) ;} \ {Z_{0, 2}^{\left(h_1\right)}  \left(y\right) = \varphi _{2}^{''}  \left(y\right) .} \end{array}
\]

Легко доказать, что верно и обратное. Другими словами, если функция $u \in  W_{p}^{\left(2, 2\right)}  \left(G\right) $ является решением задачи (1), \eqref{GrindEQ__4_}, то она является также решением задачи (1), \eqref{GrindEQ__2_} для следующих функций:
\begin{equation} \label{GrindEQ__5_}
 {\varphi _{1} \left(y\right) = Z_{0, 0}  + y Z_{0, 1}  + \int\limits_{0}^{y}  \left(y-\tau \right) Z_{0, 2}  \left(\tau \right) d\tau ;}
\end{equation}
\begin{equation} \label{GrindEQ__6_}
\varphi _{2} \left(y\right) = Z_{0, 0}^{\left(h_{ 1} \right)}  + y Z_{0, 1}^{\left(h 1\right)}  + \int\limits_{0}^{y}  \left(y-\xi \right) Z_{0, 2}^{\left(h_{ 1} \right)}  \left(\xi \right) d\xi ;
\end{equation}
\begin{equation} \label{GrindEQ__7_}
{\psi _{1} \left(x\right) = Z_{0, 0}  + x Z_{1, 0}  + \int\limits_{0}^{x}  \left(x-\eta \right) Z_{2, 0}  \left(\eta \right) d\eta ;}
\end{equation}
\begin{equation} \label{GrindEQ__8_}
 {\psi _{2} \left(x\right) = Z_{0, 0}^{\left(h_{ 2} \right)}  + x Z_{{1, 0} }^{\left(h_{ 2} \right)}  + \int\limits_{0}^{x}  \left(x-\nu \right) Z_{2, 0}^{\left(h _{2} \right)}  \left(\nu \right) d\nu}.
\end{equation}

Отметим, что функции \eqref{GrindEQ__5_} -- \eqref{GrindEQ__8_} обладают одним важным свойством, а именно, для них выполняются автоматическим образом условия согласования \eqref{GrindEQ__3_} при всех $Z_{i, j} $ и $Z_{i, j}^{(h_{k} )} $, обладающих вышеуказанными свойствами. Поэтому равенства \eqref{GrindEQ__5_} -- \eqref{GrindEQ__8_} можно рассматривать как общий вид всех функций $\varphi _{j } \left(y\right),  \psi _{j}  \left(x\right), j=\overline{1.2},$ удовлетворяющих условиям согласования   \eqref{GrindEQ__3_}.

Итак, задачи Дирихле классического вида (1), \eqref{GrindEQ__2_} и вида (1), \eqref{GrindEQ__4_} в общем случае эквивалентны. Однако задача Дирихле в неклассической трактовке (1), \eqref{GrindEQ__4_} по постановке более естественна, чем задача Дирихле (1), \eqref{GrindEQ__2_}. Это связано с тем, что в постановке  задачи Дирихле (1), \eqref{GrindEQ__4_} правые части краевых условий не требуют никаких дополнительных условий типа согласования.

Для изучения краевой задачи (1), \eqref{GrindEQ__4_} будем использовать тот факт, что любая функция $u\left(x, y\right) \in  W_{p}^{\left(2, 2\right)}  \left(G\right)$ единственным образом представима в виде:
$$
{u\left(x, y\right) = u \left(0, 0\right) + x u_{x}  \left(0, 0\right) +y u_{y}  \left(0, 0\right) + xy u_{xy } \left(0, 0\right) +}
$$
 $$
 {+\int\limits _{0}^{x}  \left(x-\alpha \right) u_{xx } \left(\alpha , 0\right) d\alpha   + y\int\limits _{0}^{x}  \left(x-\alpha \right) u_{xxy } \left(\alpha , 0\right) d\alpha  +}
 $$
\begin{equation} \label{GrindEQ__9_}
\begin{array}{c}
 {+ \int\limits _{0}^{y}  \left(y-\beta \right) u_{yy}  \left(0, \beta \right) d\beta  +x \int\limits _{0}^{y}  \left(y-\beta \right) u_{xyy}  \left(0,\beta \right) d\beta +} \\ \\\ {+ \int\limits _{0}^{x}  \int\limits _{0}^{y}  \left(x-\alpha \right) \left(y-\beta \right) u_{xxyy}  \left(\alpha ,\beta \right) d\alpha d\beta } \end{array}
\end{equation}
посредством следов $u (0,0)$, $u_{x }  \left(0, 0\right)$, $u_{y }  \left(0, 0\right)$, $u_{xy}  \left(0, 0\right)$, $u_{xx}  \left(x, 0\right), u_{xxy}  \left(x, 0\right),$\linebreak $u_{yy}  \left(0, y\right),$ $u_{xyy } \left(0, y\right)$ и старшей производной $u_{xxyy }  \left(x, y\right).$

\

{\bf 2. Операторный вид задачи Дирихле

\ \ \ \ в неклассической трактовке (1), (4)}

\

Задачу (1), (4) мы будем исследовать методом операторных
уравнений. Предварительно задачу (1), (4) запишем в виде операторного
уравнения
\begin{equation*}
Vu=Z,
\end{equation*}
где $V$ есть векторный оператор, определяемый посредством
равенства
$$
V=\left( V_{0,0} ,V_{1,0} ,V_{0,1} ,V_{2,0}
,V_{0,2} ,V_{0,0}^{(h_1)} ,V_{0,1}^{(h_1)} ,V_{0,2}^{(h_1)}, V_{0,0}^{(h_2)},V_{1,0}^{(h_2)},V_{2,0}^{(h_2)}\right):
$$$$
\qquad\qquad\qquad\qquad\qquad\qquad\qquad\qquad\qquad\qquad\qquad\qquad :W_{p}^{\left(
2,2\right) } \left( G\right) \rightarrow E_{p}^{\left(2,2\right) }
$$
а $Z$  есть заданный векторный элемент вида
$$
Z =\left(
Z _{0,0} ,Z _{1,0} ,Z _{0,1} ,Z _{2,0}
,Z _{0,2} ,Z _{0,0}^{(h_1)} ,Z _{0,1}^{(h_1)} ,Z _{0,2}^{(h_1)},Z_{0,0}^{(h_2)},Z_{1,0}^{(h_2)}, Z_{2,0}^{(h_2)}
\right)
$$
из пространства
\begin{equation*}
E_{p}^{\left(2,2\right) } \equiv R\times R\times R\times L_{p}
\left( G_{1} \right) \times L_{p} \left( G_{2} \right) \times R\times R\times L_{p}
\left( G_{2} \right) \times R\times R\times L_{p} \left( G_{1} \right) .
\end{equation*}

Заметим, что в пространстве $E_{p}^{\left(2,2\right) }$  норму
будем определять естественным образом, при помощи равенства
\begin{equation*}
\left\| Z \right\| _{E_{p}^{\left( 2,2\right) } } =\left\|
Z_{0,0} \right\| _{R} +\left\|Z_{1,0} \right\| _{R}
+\left\|Z_{0,1} \right\| _{R} +\left\| Z_{2,0}
\right\| _{L_{p} \left( G_{1} \right) } +\left\|Z_{0,2}
\right\| _{L_{p} \left( G_{2} \right) } +
\end{equation*}
\begin{equation*}
+\left\| Z_{0,0}^{(h_1)} \right\| _{R}
+\left\| Z _{0,1}^{(h_1)} \right\| _{R}
+\left\| Z_{0,2}^{(h_1)} \right\| _{L_{p} \left( G_2\right) }
+\left\| Z_{0,0}^{(h_2)} \right\| _{R} +\left\| Z_{1,0}^{(h_2)} \right\| _{R}
+\left\| Z_{2,0}^{(h_2)} \right\| _{L_{p} \left( G_1\right) }.
\end{equation*}

{\bf Определение.}\ {\it Если задача (1), (4) для любого
$Z\in E_{p}^{\left( 2,2\right) }$
имеет единственное решение $u\in W_{p}^{\left(2,2\right) } \left(
G\right)$  такое, что
\begin{equation*}
\left\| u\right\| _{W_{p}^{\left(2,2\right) } \left( G\right) }
\leq M_{1} \left\|Z\right\| _{E_{p}^{\left( 2,2\right) } } ,
\end{equation*}
то будем говорить, что оператор $V$  задачи (1),
(4) является гомеоморфизмом из
$W_{p}^{\left( 2,2\right) } \left( G\right)$  на $E_{p}^{\left(
2,2\right) }$  или задача (1), (4) везде корректно  разрешима.
Здесь  $M_{1} $- постоянное, не зависящее от $Z$.}

Очевидно,  что если оператор $V$  задачи (1), (4) является
гомеоморфизмом из $W_{p}^{\left(2,2\right) } \left( G\right)$ на
$E_{p}^{\left(2,2\right) }$, то существует ограниченный обратный
оператор
\begin{equation*}
V^{-1} :E_{p}^{\left( 2,2\right) } \rightarrow W_{p}^{\left(
2,2\right) } \left( G\right) .
\end{equation*}

\

{\bf 3. Эквивалентная система интегральных уравнений

 \ \ \ \  при исследовании задачи Дирихле (1), (4)}

\

Задачу (1), \eqref{GrindEQ__4_} мы будем изучать при помощи интегрального представления \eqref{GrindEQ__9_}   функций $u\left(x, y\right) \in W_{p}^{\left(2, 2\right)}  \left(G\right).$ Формула \eqref{GrindEQ__9_} показывает, что функция $u\left(x, y\right) \in  W_{p}^{\left(2, 2\right)}  \left(G\right)$ удовлетворяющая условиям $\left(4\right)_{1} ,$$\left(4\right)_{2, } \left(4\right)_{3, } \left(4\right)_{4} , $$\left(4\right)_{5}$ имеет следующий вид:
\begin{equation} \label{GrindEQ__10_}
\begin{array}{l} {u \left(x, y\right) = B_{0 }  \left(x, y\right) + xy b_{1, 1}  +y \int\limits _{0}^{x}  \left(x-\alpha \right) b_{2, 1}  \left(\alpha \right) d\alpha  +} \\ {+ x \int\limits _{0}^{y}  \left(y-\beta \right) b_{1, 2}  \left(\beta \right) d\beta  + \int\limits _{0}^{x}  \int\limits _{0}^{y}  \left(x-\alpha \right) \left(y-\beta \right) b_{2, 2 } \left(\alpha , \beta \right) d\alpha d\beta ,} \end{array}
\end{equation}
где
\[
B_{0}  \left(x, y\right) = Z_{0, 0}  + x Z_{1, 0}  + y Z_{0, 1 }  +\int\limits _{0}^{x}  \left(x-\alpha \right) Z_{2, 0}  \left(\alpha \right) d\alpha  +\int\limits _{0}^{y}  \left(y-\beta \right) Z_{0, 2}  \left(\beta \right) d\beta .
\]
Здесь $b_{1, 1}  \in \mathbb{R} ,  b_{2, 1}  \left(x\right) \in  L_{p}  \left(G_{1} \right), b_{1,2}  \left(y\right) \in  L_{p}  \left(G_{2} \right) $ и $b_{2, 2}  \left(x, y\right)$ неизвестные функции. Теперь функций $b_{1, 1} , b_{2, 1}  \left(x\right), b_{1, 2}  \left(y\right)$ и $b_{2,2}  \left(x, y\right)$ выберем таким образом чтобы выполнялись также условия $\left(4\right)_{6}  - \left(4\right)_{11} .$ Для этого вычислим производные $u_{x } \left(x, y\right), u_{xx } \left(x, y\right), u_{y } \left(x, y\right) $ и $u_{yy} \left(x, y\right)$ функции \eqref{GrindEQ__10_}:
\[
\begin{array}{l} {} \\ {u_{x } \left(x, y\right)=\frac{\partial  B_{0 }  \left(x, y\right)}{\partial  x}  + y b_{1, 1}  + y \int\limits _{0}^{x}  b_{2, 1}  \left(\alpha  \right) d \alpha +\int\limits _{0}^{y}  \left(y-\beta \right) b_{1, 2}  \left(\beta \right) d \beta  + {\rm \; \; \; \; \; }} \\ {} \end{array}
\]
$$
+\int\limits _{0}^{x}\int\limits _{0}^{y} \left(y-\beta \right) b_{2, 2}  \left(\alpha ,\beta \right) d\alpha  d \beta ,\ \ \mbox{где}\ \
\frac{\partial  B_{0 }  \left(x, y\right)}{\partial  x} = Z_{1, 0}  + \int\limits _{0}^{x}  Z_{2, 0}  \left(\alpha \right) d \alpha ;
$$
\[
u_{xx } \left(x, y\right) =
{\frac{\partial^{2}  B_{0 }  \left(x, y\right)}{\partial  x_{}^{2} }  + y b_{2, 1}  \left(x\right) }
+
\]\[
+\int\limits _{0}^{y}  \left(y-\beta \right) b_{2, 2}  \left(x,  \beta \right) d \beta , \ \ \mbox{где}\ \
\frac{\partial^{2}  B_{0} \left(x, y\right)}{\partial  x^{2} } =Z_{2, 0}  \left(x\right);
\]
\[
u_{y}  \left(x,y\right)= \frac{\partial  B_{0 } \left(x,y\right)}{\partial  y} +xb_{1, 1} +\int\limits _{0}^{x}  \left(x-\alpha \right) b_{2, 1}  \left(\alpha \right) d \alpha +
x \int\limits _{0}^{y}  b_{1, 2}  \left(\beta \right) d \beta +
 \]
\[
 {+\int\limits _{0}^{x} \int\limits _{0}^{y}  \left(x-\alpha \right) b_{2, 2}  \left(\alpha , \beta \right) d \alpha  d \beta , }
\ \
\mbox{где}
\ \
 \frac{\partial  B_{0}  \left(x, y\right)}{\partial  y} =Z_{0, 1 } +\int\limits _{0}^{y}  Z_{0, 2}  \left(\beta \right) d \beta ;
 \]
 \[
 u_{yy } \left(x, y\right)=\frac{\partial^{2}  B_{0 }  \left(x, y\right)}{\partial  y^{2} }  +
 {xb_{1, 2}  \left(y\right)+ \int\limits _{0}^{x}  \left(x-\alpha \right) b_{2, 2} \left(\alpha , y\right) d \alpha , }
\]
где
\[
\frac{\partial _{}^{2}  B_{0 }  \left(x, y\right)}{\partial  y_{}^{2} } = Z_{0, 2}  \left(y\right).
\]
Выражения этих производных показывает что, справедлива следующие равенства
\[
{u \left(h_{1} , 0\right)=B_{0}  \left(h_{1},  0\right) = Z_{0, 0} + h_{1}    Z_{1, 0 } +}
+\int\limits _{0}^{h_{ 1} } \left(h_{1} -\alpha \right) Z_{2, 0}  \left(\alpha \right) d\alpha ,
\]
\[
{u_{y} \left(h_{1} , 0\right) = Z_{0, 1} + b_{1, 1 }  h_{1}  +}
{\int\limits _{0}^{h_{ 1} } \left(h_{1} -\alpha \right)b_{2, 1 } \left(\alpha \right) d \alpha , }
\]
\[
u_{yy} \left(h_{1} , y\right) = Z_{0, 2} \left(y\right) +
{h_{1}    b_{1, 2 } \left(y\right) +\int\limits _{0}^{h_{ 1} } \left(h_{1} -\alpha \right) b_{2, 2}  \left(\alpha , y\right) d \alpha , }
\]
\[u \left(0,h_{2} \right) = {B_{0} \left(0, h_{2} \right)= Z_{0, 0}  +h_{2}   Z_{0, 1} + \int\limits _{0}^{h_{ 2} } \left(h_{2} - \beta \right) Z_{0, 2}  \left(\beta \right) d\beta ,}
\]
\[ {u_{x} \left(0, h_{2} \right)=Z_{1, 0} +h_{2} b_{1,1} +\int\limits _{0}^{h_{ 2} }  \left(h_{2} -\beta \right)b_{1, 2} \left(\beta \right) d \beta ,}
\]
\[
{u_{xx} \left(x, h_{2} \right)=Z_{2, 0} \left(x\right)+h_{2}  b_{2,1} \left(x\right)+\int\limits _{0}^{h _{2} } \left(h_{2} -\beta \right) b_{2, 2} \left(x, \beta \right) d\beta .}
\]

Поэтому условия \eqref{GrindEQ__4_}$_{6} $- \eqref{GrindEQ__4_}$_{11} $ можно записать в следующем виде:
\[
\begin{array}{l} {Z_{0, 0} +h_{1}  Z_{1, 0} +\int\limits _{0}^{h_{ 1} } \left(h_{1} -\alpha \right) Z_{2, 0} \left(\alpha \right) d \alpha =Z_{0, 0}^{\left(h _{1} \right)} ;} \\ {Z_{0, 1} +h_{1}  b_{1, 1} +\int\limits _{0}^{h_{ 1} } \left(h_{1} -\alpha \right) b_{2, 1} \left(\alpha \right)d \alpha =Z_{0, 1}^{\left(h _{1} \right)} ;} \\ {Z_{0, 2 } \left(y\right)+h_{1}  b_{1, 2} \left(y\right)+\int\limits _{0}^{h_{ 1} } \left(h_{1} -\alpha \right)b_{2, 2 } \left(\alpha , y\right) d \alpha =Z_{0, 2}^{\left(h _{1} \right)} \left(y\right);} \\ {Z_{0, 0 } +h_{2}  Z_{0, 1} +\int\limits _{0}^{h_{ 2} } \left(h_{2} -\beta \right)Z_{0, 2} \left(\beta \right) d \beta =Z_{0, 0}^{\left(h _{2} \right)} ;} \\ {Z_{1,0 } +h_{2}  b_{1, 1} +\int\limits _{0}^{h_{ 2} } \left(h_{2} -\beta \right) b_{1, 2 } \left(\beta \right) d \beta =Z_{1, 0}^{\left(h _{2} \right)} ;} \\ {Z_{2, 0 } \left(x\right)+h_{2}  b_{2, 1} \left(x\right)+\int\limits _{0}^{h_{ 2} } \left(h_{2} -\beta \right)b_{2, 2 } \left(x, \beta \right) d \beta =Z_{2, 0}^{\left(h _{2} \right)} \left(x\right).} \end{array}
\]

Таким образом, мы получили следующая система интегральных уравнений:
\begin{equation} \label{GrindEQ__11_}
\left\{\begin{array}{l} {Z_{0, 1} +h_{1}  b_{1, 1} +\int\limits _{0}^{h_{ 1} } \left(h_{1} -\alpha \right) b_{2, 1} \left(\alpha \right)d \alpha =Z_{0, 1}^{\left(h _{1} \right)} ,} \\ {Z_{0, 2 } \left(y\right)+h_{1}  b_{1, 2} \left(y\right)+\int\limits _{0}^{h_{ 1} } \left(h_{1} -\alpha \right)b_{2, 2 } \left(\alpha , y\right) d \alpha =Z_{0, 2}^{\left(h _{1} \right)} \left(y\right),} \\ {Z_{1,0 } +h_{2}  b_{1, 1} +\int\limits _{0}^{h_{ 2} } \left(h_{2} -\beta \right) b_{1, 2 } \left(\beta \right) d \beta =Z_{1, 0}^{\left(h _{2} \right)} ,} \\ {Z_{2, 0 } \left(x\right)+h_{2}  b_{2, 1} \left(x\right)+\int\limits _{0}^{h_{ 2} } \left(h_{2} -\beta \right)b_{2, 2 } \left(x, \beta \right) d \beta =Z_{2, 0}^{\left(h _{2} \right)} \left(x\right).} \end{array}\right.
\end{equation}

Здесь $b_{1, 1} $, $b_{1, 2} \left(y\right)$, $b_{2, 1} \left(x\right)$ и $b_{2, 2} \left(x,y\right)$- искомые элементы.

Теперь потребуем чтобы, функция \eqref{GrindEQ__10_} являлась решением уравнения (1). Для этого выражение \eqref{GrindEQ__10_} этой функции учтем в уравнении (1).

Тогда после замены
\[
u\left(x, y\right)=B_{0} \left(x, y\right)+\widetilde{u}\left(x, y\right),
\]
где
\[
{\widetilde{u}\left(x, y\right)=x y b_{1, 1} +y \int\limits _{0}^{x}  \left(x-\alpha \right) b_{2, 1} \left(\alpha \right) d \alpha +}
 { x\int\limits _{0}^{y} \left(y-\beta \right) b_{1, 2}  \left(\beta  \right) d \beta +}
\]
\[
+{\int\limits _{0}^{x} \int\limits _{0}^{y}  \left(x-\alpha \right) \left(y-\beta \right) b_{2, 2}  \left(\alpha , \beta \right) d\alpha d\beta }
\]
уравнение (1) можно записать в виде
\begin{equation} \label{GrindEQ__12_}
\left(V_{2, 2}  \widetilde{u}\right) \left(x, y\right)=\widetilde{R} \left(x, y\right),
\end{equation}
где
\[
\widetilde{R} \left(x, y\right)=Z_{2, 2 } \left(x, y\right)-\left(V_{2, 2}  B_{0} \right) \left(x, y\right) .
\]

Очевидно, что производные функции $\widetilde{u}$ можно вычислить посредством равенств
\[
\widetilde{u}_{x} \left(x, y\right)=y b_{1, 1} +y \int\limits _{0}^{x}  b_{2, 1} \left(\alpha \right) d \alpha +
\]
\[
+ \int\limits _{0}^{y} \left(y-\beta \right) b_{1, 2}  \left(\beta  \right) d \beta +
\int\limits _{0}^{x} \int\limits _{0}^{y}  \left(y-\beta \right) b_{2, 2}  \left(\alpha , \beta \right) d\alpha d\beta ,
\]
\[
 \widetilde{u}_{y} \left(x, y\right)=x b_{1, 1} + \int\limits _{0}^{x}  \left(x-\alpha \right) b_{2, 1} \left(\alpha \right) d \alpha
 +
 \]
 \[
 + x\int\limits _{0}^{y} b_{1, 2}  \left(\beta  \right) d \beta +
 \int\limits _{0}^{x} \int\limits _{0}^{y}  \left(x-\alpha \right) b_{2, 2}  \left(\alpha , \beta \right) d\alpha d\beta ,
 \]
 \[
 \widetilde{u}_{xx} \left(x, y\right)=y b_{2, 1} \left(x\right)+ \int\limits _{0}^{y}  (y-\beta )b_{2, 2} \left(x, \beta \right) d\beta ,
 \]
 \[
 \widetilde{u}_{yy} \left(x, y\right)=x b_{1, 2} \left(y\right)+ \int\limits _{0}^{x}  (x-\alpha )b_{2, 2} \left(\alpha ,y \right) d\alpha ,
 \]
 \[
 \widetilde{u}_{xy} \left(x, y\right)=b_{1, 1} + \int\limits _{0}^{x}  b_{2, 1} \left(\alpha \right) d \alpha + \int\limits _{0}^{y} b_{1, 2}
  \left(\beta \right) d \beta +
\int\limits _{0}^{x} \int\limits _{0}^{y}  b_{2, 2}  \left(\alpha , \beta \right) d\alpha d\beta ,
\]
\[
{\widetilde{u}_{xxy} \left(x, y\right)=b_{2, 1} \left(x\right)+ \int\limits _{0}^{y}  b_{2, 2} \left(x, \beta \right) d\beta ,}
\]
\[ {\widetilde{u}_{xyy} \left(x, y\right)= b_{1, 2} \left(y\right)+ \int\limits _{0}^{x} b_{2, 2} \left(\alpha ,y \right) d\alpha ,}
\
\  {\widetilde{u}_{xxyy} \left(x, y\right)=b_{2, 2} \left(x, y\right)}.
 \]

Поэтому уравнение \eqref{GrindEQ__12_} можно привести к виду
\[
\left(Nh\right)\left(x,y\right)\equiv b_{2, 2} \left(x,y\right)+x y b_{1, 1}  a_{0, 0 } \left(x, y\right)+y\int\limits _{0}^{x}  \left(x-\alpha \right) a_{_{0, 0} } \left(x, y\right) b_{2, 1} \left(\alpha \right) d\alpha +
\]
\[
+ x\int\limits _{0}^{y}  \left(y-\beta \right) a_{_{0, 0} } \left(x, y\right) b_{1, 2} \left(\beta \right) d\beta +
 \]
\[
{+\int\limits _{0}^{x}  \int\limits _{0}^{y}  \left(x-\alpha \right)\left(y-\beta \right) a_{0, 0 } \left(x, y\right)b_{2, 2} \left(\alpha ,\beta \right)d\alpha d\beta +yb_{1, 1}^{}  a_{1,0}^{}  \left(x, y\right)+}
\]
\[
+{y\int\limits _{0}^{x}  a_{1, 0 } \left(x, y\right)b_{2, 1} \left(\alpha \right)d\alpha +\int\limits _{0}^{y}  \left(y-\beta \right) a_{_{1, 0} } \left(x, y\right) b_{1, 2} \left(\beta \right) d\beta +}
\]
\[ {+\int\limits _{0}^{x}  \int\limits _{0}^{y}  \left(y-\beta \right) a_{1, 0 } \left(x, y\right)b_{2, 2} \left(\alpha ,\beta \right)d\alpha d\beta +}
 x b_{1, 1}   a_{0, 1 } \left(x, y\right)+
  \]
 \[
 +\int\limits _{0}^{x}  \left(x-\alpha \right) a_{_{0, 1} } \left(x, y\right)b_{2, 1 } \left(\alpha \right) d\alpha +x\int\limits _{0}^{y} a_{_{0, 1} } \left(x, y\right) b_{1, 2} \left(\beta \right) d\beta +
 \]
 \[ {+\int\limits _{0}^{x}  \int\limits _{0}^{y}  \left(x-\alpha \right) a_{0, 1 } \left(x, y\right)b_{2, 2} \left(\alpha ,\beta \right)d\alpha d\beta }
  +y{a_{2, 0 } \left(x, y\right)  b_{2, 1} (x)+}
\]
\[+ \int\limits _{0}^{y}  \left(y-\beta \right) a_{2, 0 } \left(x, y\right)b_{2, 2} \left(x,\beta \right)d\beta +
x  a_{0, 2 } \left(x, y\right)  b_{1, 2} \left(y\right)+
\]
\[
+\int\limits _{0}^{x}  \left(x-\alpha \right) a_{0, 2 } \left(x, y\right)b_{2, 2} \left(\alpha ,y\right)d\alpha +  b_{1, 1}  a_{1, 1 } \left(x, y\right)+\int\limits _{0}^{x}  a_{_{1, 1} } \left(x, y\right)b_{2, 1 } \left(\alpha \right) d\alpha +
 \]
 \[
 {+ \int\limits _{0}^{y}  a_{_{1, 1} } \left(x, y\right)b_{1, 2 } \left(\beta \right) d\beta } {+\int\limits _{0}^{x}  \int\limits _{0}^{y}  a_{1, 1 } \left(x, y\right)b_{2, 2} \left(\alpha ,\beta \right)d\alpha d\beta +}
 \]
 \[
 {+a_{2, 1 } \left(x, y\right)b_{2, 1} (x)}  {+\int\limits _{0}^{y}  a_{_{2, 1} } \left(x, y\right)b_{2, 2 } \left(x,\beta \right) d\beta }  {+a_{_{1, 2} } \left(x, y\right)b_{1, 2 } (y)+}
 \]
 \[
  {+\int\limits _{0}^{x} a_{_{1, 2} } \left(x, y\right)b_{2, 2 } \left(\alpha , y\right) d\alpha =\widetilde{R}\left(x, y\right),}
 \]
где
\[
h=\left(b_{1, 1 ,  }  b_{2,1} \left(x\right),   b_{1, 2}  \left(y\right), b_{2,2} \left(x, y\right)\right) \in  \widetilde{E}_{p}^{\left(2, 2\right)} \equiv R\times L_{p} \left(G_{1} \right)\times L_{p} \left(G_{2} \right)\times L_{p} \left(G\right)
\]
неизвестная четверка.

Произведя здесь некоторые группировки, имеем
\[
\left(Nh\right)\left(x,y\right) =\int\limits _{0}^{x} \left[\left(x-\alpha  \right) \left(y a_{0, 0 } \left(x,y\right)+a_{0, 1} \left(x,y\right)\right)+
\right.\]
\[\left. +y a_{1, 0} \left(x, y\right)+a_{1, 1} \left(x,y\right)\right] b_{2, 1}  \left(\alpha \right) d\alpha +
\]
\[
{+\int\limits _{0}^{x} [\left(x-\alpha  \right) a_{0, 2 } \left(x,y\right)+a_{1, 2} \left(x,y\right)] b_{2, 2 } \left(\alpha , y\right)d \alpha +}
\]
\[ {+\int\limits _{0}^{y} [\left(y-\beta  \right) \left(x a_{0, 0} \left(x,y\right)+a_{1, 0} \left(x,y\right)\right)}  {+x a_{0, 1} \left(x, y\right)+a_{1, 1} \left(x,y\right)] b_{1, 2}  \left(\beta \right) d\beta +}
\]
\[ {+\int\limits _{0}^{y} [\left(y-\beta  \right) a_{2, 0 } \left(x,y\right)+a_{2, 1} \left(x,y\right)] b_{2, 2 } \left(x, \beta \right)d \beta +}
\]
\[ {+\int\limits _{0}^{x}  \int\limits _{0}^{y}  [\left(x-\alpha \right)\left(y-\beta \right) a_{0, 0 } \left(x, y\right)+\left(y-\beta \right) a_{1,0}^{}  \left(x, y\right)+}
\]
\[ {+\left(x-\alpha \right) a_{_{0, 1} } \left(x, y\right)+a_{1, 1} \left(x,y\right)] b_{2, 2} \left(\alpha , \beta \right) d\alpha d\beta +}
\]
\[ {+[x y a_{0, 0} \left(x,y\right)+y a_{1, 0} \left(x,y\right)+x a_{0, 1} \left(x,y\right)+}a_{1, 1} \left(x,y\right)] b_{1, 1} +
\]
\[
+ {[y a_{2, 0} \left(x,y\right)+a_{2, 1} \left(x,y\right)] b_{2, 1 } \left(x\right)+}
\]
\begin{equation} \label{GrindEQ__13_}
{+[x a_{0,2} \left(x,y\right)+a_{1, 2} \left(x,y\right)] b_{1, 2} \left(y\right)} {+b_{2, 2} \left(x, y\right)=\widetilde{R}\left(x, y\right)}.
\end{equation}

Таким образом, решение задачи (1), \eqref{GrindEQ__4_} в пространстве $W_{p}^{\left(2, 2\right)} \left(G\right)$ сведено к решению систем интегральных уравнений \eqref{GrindEQ__11_}, \eqref{GrindEQ__13_} относительно четверки неизвестных

$h=\left(b_{1, 1} , b_{2, 1} \left(x\right), b_{1, 2} \left(y\right)b_{2, 2} \left(x, y\right)\right)$ в пространстве $\widetilde{E}_{p}^{(2, 2)} $.

Система интегральных уравнений \eqref{GrindEQ__11_} можно разрешить относительно $b_{2, 1} \left(x\right),$   $b_{1, 2} \left(y\right)$,    $ b_{1, 1} $     следующим образом:
\begin{equation} \label{GrindEQ__14_}
{b_{2, 1} \left(x\right)=\frac{Z_{2, 0}^{\left(h_{2} \right)}  \left(x\right)-Z_{2, 0} \left(x\right)}{h_{2} }}  {- \frac{1}{h_{2} } \int\limits _{0}^{h_{ 2} } \left(h_{ 2} -\beta \right) b_{2, 2} \left(x, \beta \right) d\beta ,}
\end{equation}
\begin{equation} \label{GrindEQ__15_}
 {b_{1, 2} \left(y\right)=\frac{Z_{0, 2}^{\left(h_{ 1} \right)}  \left(y\right)-Z_{0, 2} \left(y\right)}{h_{1} } - }  { \frac{1}{h_{1} } \int\limits _{0}^{h_{ 1} } \left(h_{1} -\alpha \right) b_{2, 2} \left(\alpha , y\right) d\alpha ,}
\end{equation}
\[
b_{1, 1} =\frac{Z_{0, 1}^{\left(h_{ 1} \right)} -Z_{0, 1} }{h_{1} } - \frac{1}{h_{1} } \int\limits _{0}^{h_{ 1} } \left(h_{1} -\alpha \right) b_{2, 1} \left(\alpha \right) d\alpha ,
\]
или же,
\begin{equation} \label{GrindEQ__16_}
b_{1, 1} =\frac{Z_{1, 0}^{\left(h_{ 2} \right)} -Z_{1, 0} }{h_{2} } - \frac{1}{h_{2} } \int\limits _{0}^{h_{ 2} } \left(h_{2} -\beta \right) b_{1, 2} \left(\beta \right) d\beta ,
\end{equation}

Таким образом, при вышеналоженных условиях доказана следующая

{\bf Теорема 1.} {\it Для везде корректно разрешимости задачи (1), (4) в пространстве $W_{p}^{\left(2, 2\right)} \left(G\right)$ необходимо и достаточно чтобы система интегральных уравнений
\eqref{GrindEQ__13_}, \eqref{GrindEQ__14_}, \eqref{GrindEQ__15_}, \eqref{GrindEQ__16_} была везде корректно разрешимой в пространстве $\widetilde{E}_{p}^{(2, 2)} $.}\\

{\bf 4. Построение эквивалентного интегрального уравнения}

\

Из \eqref{GrindEQ__14_} -- \eqref{GrindEQ__15_} следует, что  $b_{1, 1} $ можно представить в виде
\begin{equation} \label{GrindEQ__17_}
b_{1, 1} =\frac{Z_{0, 1}^{\left(h_{ 1} \right)} -Z_{0, 1} }{h_{1} }  -
 $$$$
 - \frac{1}{h_{1} } \int\limits _{0}^{h_{ 1} } \left(h_{1} -\alpha \right) \left[\frac{Z_{2, 0}^{\left(h_{ 2} \right)} \left(\alpha \right)-Z_{2, 0} \left(\alpha \right)}{h_{ 2} } -\frac{1}{h_{2} } \int\limits _{0}^{h_{ 2} } \left(h_{2} -\beta \right) b_{2, 2}  \left(\alpha , \beta \right) d\beta \right] d\alpha ,
\end{equation}
или же
\[
b_{1, 1} =\frac{Z_{1, 0}^{\left(h_{ 2} \right)} -Z_{1, 0} }{h_{2} } -
\]
$$
- \frac{1}{h_{2} } \int\limits _{0}^{h_{ 2} } \left(h_{2} -\beta \right) \left[ \frac{Z_{0, 2}^{\left(h_{ 1} \right)} \left(\beta \right)-Z_{0, 2} \left(\beta \right)}{h_{ 1} } -\frac{1}{h_{1} } \int\limits _{0}^{h_{ 1} } \left(h_{1} -\alpha \right) b_{2, 2}  \left(\alpha , \beta \right) d\alpha \right] d\beta  .
$$

Тогда очевидно, что учитывая выражении \eqref{GrindEQ__14_}, \eqref{GrindEQ__15_} и \eqref{GrindEQ__17_} неизвестных       $b_{2, 1} \left(x\right)$, $b_{1, 2} \left(y\right)$ и $b_{1, 1} $ уравнение \eqref{GrindEQ__13_} можно привести также к самостоятельному интегральному уравнению с неизвестной функцией $b_{2, 2} \left(x, y\right)$:
\[
\left(Nb_{2, 2} \right)\left(x,y\right)\equiv
\]
\[
\equiv\int\limits _{0}^{x} \left[\left(x-\alpha  \right) \left(y a_{0, 0 } \left(x,y\right)+a_{0, 1} \left(x,y\right)\right)+y a_{1, 0} \left(x, y\right)+a_{1, 1} \left(x,y\right)\right] \times
\]
\[
\times\left[\frac{Z_{2, 0}^{\left(h_{ 2} \right)} \left(\alpha \right)-Z_{2, 0} \left(\alpha \right)}{h_{ 2} }
-\frac{1}{h_{2} } \int\limits _{0}^{h_{ 2} } \left(h_{2} -\beta \right) b_{2, 2}  \left(\alpha , \beta \right) d\beta \right] d\alpha  +
\]
\[
{+\int\limits _{0}^{x} \left[\left(x-\alpha  \right) a_{0, 2 } \left(x,y\right)+a_{1, 2} \left(x,y\right)\right] b_{2, 2 } \left(\alpha , y\right)d \alpha +}
\]
\[
+\int\limits _{0}^{y} [\left(y-\beta  \right) \left(x a_{0, 0} \left(x,y\right)+a_{1, 0} \left(x,y\right)\right)+
x a_{0, 1} \left(x, y\right)+a_{1, 1} \left(x,y\right)]\times
\]
\[
\times \left[ \frac{Z_{0, 2}^{\left(h_{ 1} \right)} \left(\beta \right)-Z_{0, 2} \left(\beta \right)}{h_{ 1} }
 -\frac{1}{h_{1} } \int\limits _{0}^{h_{ 1} } \left(h_{1} -\alpha \right) b_{2, 2}  \left(\alpha , \beta \right) d\alpha \right] d\beta  +
\]
\[
{+\int\limits _{0}^{y} [\left(y-\beta  \right) a_{2, 0 } \left(x,y\right)+a_{2, 1} \left(x,y\right)] b_{2, 2 } \left(x, \beta \right)d \beta +}
\]
\[
{+ \int\limits _{0}^{x}  \int\limits _{0}^{y}  [\left(x-\alpha \right)\left(y-\beta \right) a_{0, 0 } \left(x, y\right)+\left(y-\beta \right) a_{1,0}^{}  \left(x, y\right)+}
\]
\[
 {+\left(x-\alpha \right) a_{_{0, 1} } \left(x, y\right)+a_{1, 1} \left(x,y\right)] b_{2, 2} \left(\alpha , \beta \right) d\alpha d\beta +}
 \]
 \[ +[x y a_{0, 0} \left(x,y\right)+y a_{1, 0} \left(x,y\right)+x a_{0, 1} \left(x,y\right)+a_{1, 1} \left(x,y\right)] \times
 \]
 \[
 \times  \left[\frac{Z_{0, 1}^{\left(h_{ 1} \right)} -Z_{0, 1} }{h_{1} } - \frac{1}{h_{1} } \int\limits _{0}^{h_{ 1} } \left(h_{1} -\alpha \right) \left(\frac{Z^{\left(h_{2} \right)}_{2, 0} \left(\alpha \right)  -Z_{2, 0} \left(\alpha \right)}{h_{2} } -\right.\right.
\]
\[
\left. \left.- \frac{1}{h_{2} } \int\limits _{0}^{h_{ 2} } \left(h_{ 2} -\beta \right) b_{2, 2} \left(\alpha , \beta \right) d\beta \right) d\alpha \right]+[y a_{2, 0} \left(x,y\right)+a_{2, 1} \left(x,y\right)] \times
\]
\[ \times  \left[\frac{Z_{2, 0}^{\left(h_{2} \right)}  \left(x\right)-Z_{2, 0} \left(x\right)}{h_{2} } -  \frac{1}{h_{2} } \int\limits _{0}^{h_{ 2} } \left(h_{ 2} -\beta \right) b_{2, 2} \left(x, \beta \right) d\beta \right] +
\]
\[
+ \left[x a_{0,2} \left(x,y\right)+a_{1, 2} \left(x,y\right)\right] \left[\frac{Z_{0, 2}^{\left(h_{ 1} \right)}  \left(y\right)-Z_{0, 2} \left(y\right)}{h_{1} } -\right.
\]
\[
\left.- \frac{1}{h_{1} } \int\limits _{0}^{h_{ 1} } \left(h_{1} -\alpha \right) b_{2, 2} \left(\alpha , y\right) d\alpha \right]+b_{2,2} \left(x, y\right)=\widetilde{R}\left(x, y\right).\eqno(18)
\]

Оператор $N$ уравнения (18) линеен. Используя условия, наложенные
на коэффициенты  $ a_{i,j} $, можно доказать, что этот оператор
является ограниченным оператором из $L_{p} (G)$  в  $L_{p} (G)$,
$1\leq p\leq \infty .$

При помощи, например, метода
последовательных приближений можно доказать, что уравнение (18) для
любой правой части $\widetilde{R} \left( x,y\right) \in L_{p} \left(
G\right)$ имеет единственное решение $b_{2,2}(x,y)\in L_{p} \left( G\right) ,$
где $1\leq p\leq \infty ,$ и это решение удовлетворяет условию
\begin{equation*}
\left\| b_{2,2}\right\| _{L_{p\left( G\right) } } \leq M_{2} \left\|
\widetilde{R} \right\| _{L_{p\left( G\right) } },
\end{equation*}
где    $M_{2} $-постоянное, не зависящее от  $\widetilde{R}$. Далее,
очевидно, что если  $b_{2,2}\in L_{p}
\left( G\right)$  есть решение уравнения (18), то решение задачи
(1), (4) можно найти при помощи равенства (10) с учетом (14)--(16).

Поэтому справедлива

{\bf Теорема 2.}\ {\it Оператор $V$ задачи (1), (4) есть
гомеоморфизмом из\linebreak $W_{p}^{\left( 2,2\right) } \left( G\right)$ на
$E_{p}^{\left(2,2\right) }$.}

Отметим, что в некоторые краевые задачи в неклассических трактовках  исследованы  в работах автора [6-7].

\

{\bf 5. Выводы}

1. В работе получена эквивалентная система интегральных уравнений типа Фредгольма при исследовании задачи Дирихле для обобщенного уравнения Манжерона с негладкими коэффициентами  в неклассической трактовке (1), (4).

2. При негладких условиях на коэффициенты уравнения, в прямоугольной области для этой задачи найдены  условия корректной разрешимости в интегральном виде на основе метода интегральных представлений.

\newpage

\begin{titlepage}

\centerline{\bf Сведения об авторе}

\

{\bf Фамилия, имя, отчество:} {\bf Мамедов Ильгар Гурбат оглы }

{\bf Название организации:} Институт Кибернетики им. А.И.Гусейнова НАН Азербайджана

{\bf Должность, уч. степень, уч. звание:} ведущий научный сотрудник, кандидат физико-математических наук, доцент.

{\bf Почтовый адрес:} Азербайджан, AZ 1141, г. Баку, ул. Б. Вагабзаде, 9,

{\bf Телефон, Email, fax:} (99412) 539 28 26, ilgar-mammadov@rambler.ru, (99412) 539 28 26.
\end{titlepage}

\end{document}